\documentclass[12pt]{article}
\usepackage{amsfonts}
\usepackage{graphicx}
\usepackage{amsmath}
\usepackage[francais]{babel}

\newtheorem{theorem}{Theorem}

\newtheorem{example}[theorem]{Exemple}

\newtheorem{proposition}[theorem]{Proposition}

\newtheorem{theoreme}[theorem]{Théorème}

\newenvironment{demonstration}[1][Démonstration]{\textbf{#1:} }
{\ \rule{0.5em}{0.5em}}
\input{tcilatex}

\begin{document}

\begin{center}
{\LARGE Feuilletage canonique sur le fibr\'{e} de Weil}

{\LARGE \bigskip }

\textbf{Basile Guy Richard BOSSOTO}

\textbf{\medskip }

Universit\'{e} Marien NGOUABI

Facult\'{e} des Sciences, D\'{e}partement de Math\'{e}matiques

BP: 69, Brazzaville, Congo.

E-mail : bossotob@yahoo.fr

\bigskip
\end{center}

\textbf{\ Abstract }: Let be $M$ a smooth manifold, $A$ a local algebra and $%
M^{A}$ a manifold of infinitely near points on $M$ of kind $A$. We build the
canonical foliation on $M^{A}$ et we show that the canonical foliation on
the tangent bundle $TM$ is the foliation defined by his canonical field.

\bigskip

\textbf{R\'{e}sum\'{e} }: Soit $M$ une vari\'{e}t\'{e} diff\'{e}rentielle, $%
A $ une alg\`{e}bre locale et $M^{A}$ la vari\'{e}t\'{e} des points proches
de $M$ d'esp\`{e}ce $A$. Nous construisons le feuilletage canonique sur $%
M^{A}$ et montrons que le feuilletage canonique du fibr\'{e} tangent $TM$,
est le feuilletage d\'{e}fini par son champ canonique.

\medskip


\textbf{Mots cl\'{e}s}\textit{\ : }Points proches, alg\`{e}bre locale,
feuilletages, d\'{e}rivations.

\section{Pr\'{e}liminaires}

\subsection{D\'{e}rivation d'une alg\`{e}bre}

Soit $A$ une alg\`{e}bre commutative et unitaire sur $\mathbb{R}$ et $%
\mathfrak{M}$ un $A$-module. Une d\'{e}rivation de $A$ dans $\mathfrak{M}$
est une application $\mathbb{R}$-lin\'{e}aire $d:A\longrightarrow \mathfrak{M%
}$ telle que 
\begin{equation*}
d(ab)=d(a)b+ad(b)
\end{equation*}%
pour tous $a,b\in A$. Evidemment $d\left( \lambda \right) =0$ pour tout $%
\lambda \in 
\mathbb{R}
$.

On note $Der(A,\mathfrak{M})$ le $A$-module des d\'{e}rivations de $A$ dans $%
\mathfrak{M}$. Lorsque $\mathfrak{M}=A$, une d\'{e}rivation de $A$ dans $A$
est simplement appel\'{e}e d\'{e}rivation de $A$ et on note $Der_{\mathbb{R}%
}(A)$ o\`{u} simplement $Der(A)$ s'il n y a pas de confusion, le $A$-module
des d\'{e}rivations de $A$.

Si $A$ et $B$ sont deux alg\`{e}bres quelconque et si $\varphi
:A\longrightarrow B$ est un homomorphisme d'alg\`{e}bres, alors $B$ est un $%
A $-module.

Soit $\varphi :A\longrightarrow B$ un homomorphisme d'alg\`{e}bres. Une
application $\mathbb{R}$-lin\'{e}aire $d:A\longrightarrow B$ est une $%
\varphi $-d\'{e}rivation si

\begin{equation*}
d(ab)=d(a)\varphi(b)+\varphi(a)d(b)
\end{equation*}
pour tous $a,b\in A$.

\subsection{\textbf{\ }Alg\`{e}bre locale et vari\'{e}t\'{e} des points
proches}

Une alg\`{e}bre locale au sens de Weil est une alg\`{e}bre r\'{e}elle
commutative unitaire $A$, de dimension finie sur $\mathbb{R}$, ayant un id%
\'{e}al maximal unique $\mathfrak{m}$ de codimension $1$. On a ainsi 
\begin{equation*}
A=\mathbb{R}\oplus \mathfrak{m}\text{.}
\end{equation*}%
Dans ce cas, compte tenu du lemme de Nakayama, l'id\'{e}al maximal $%
\mathfrak{m}$ est nilpotent. Le plus petit entier positif $k$ tel que $%
\mathfrak{m}^{k+1}=(0)$ est la hauteur de $A$ et la dimension sur $\mathbb{R}
$ de $\mathfrak{m}/\mathfrak{m}^{2}$ est la largeur ou la profondeur de $A$.

Lorsque $A$ est une alg\`{e}bre locale (de dimension finie), l'ensemble $%
Der\left( A\right) $ des d\'{e}rivations de $A$ est une alg\`{e}bre de Lie
de dimension finie : c'est l'alg\`{e}bre de Lie du groupe de Lie $Aut(A)$
des automorphismes de $A$.

Dans toute la suite, $M$ est une vari\'{e}t\'{e} diff\'{e}rentielle
paracompacte de classe $C^{\infty }$ de dimension $n$ et\ $A$\ une alg\`{e}%
bre locale au sens de Weil. On note $C^{\infty }(M)$ l'alg\`{e}bre des
fonctions num\'{e}riques de classe $C^{\infty }$ sur $M$, $\mathfrak{X}(M)$
le $C^{\infty }(M)$-module des champs de vecteurs sur $M$ \ et $TM$ l'espace
tangent \`{a} $M$.

Un point proche de $p\in M$ d'esp\`{e}ce $A$\cite{wei}, est un homomorphisme
d'alg\`{e}bres 
\begin{equation*}
\xi :C^{\infty }(M)\longrightarrow A
\end{equation*}%
tel que, pour tout $f\in C^{\infty }(M)$, 
\begin{equation*}
\left[ \xi (f)-f(p)\right] \in \mathfrak{m}\text{.}
\end{equation*}

On note $M_{p}^{A}$ l'ensemble des points proches de $p\in M$ d'esp\`{e}ce $%
A $ et 
\begin{equation*}
M^{A}=\dbigcup\limits_{p\in M}M_{p}^{A}\text{.}
\end{equation*}

L'ensemble $M^{A}=Hom_{A\lg }(C^{\infty }(M),A)$ est une vari\'{e}t\'{e} diff%
\'{e}rentielle de dimension $\dim (M)\cdot \dim (A)$ et est appel\'{e} vari%
\'{e}t\'{e} des points proches de $M$ d'esp\`{e}ce $A$\cite{wei} ou
simplement fibr\'{e} de Weil d'esp\`{e}ce $A$.

\begin{enumerate}
\item Lorsque $A=%
\mathbb{R}
$, on identifie $M^{%
\mathbb{R}
}$ \`{a} $M$ par l'application 
\begin{equation*}
M\longrightarrow M^{%
\mathbb{R}
}=Hom_{A\lg }(C^{\infty }(M),%
\mathbb{R}
),p\longmapsto \left\{ f\longmapsto f(p)\right\} \text{.}
\end{equation*}

\item Lorsque $V$ est un espace vectoriel de dimension finie $p$ dont une
base est $(v_{1},...;v_{p})$, si $(v_{1}^{\ast },...,v_{p}^{\ast })$ d\'{e}%
signe la base duale de $(v_{1},...;v_{p})$, alors l'application 
\begin{equation*}
V^{A}\overset{\theta }{\longrightarrow }V\otimes A,\xi \longmapsto \underset{%
i=1}{\overset{p}{\sum }}v_{i}\otimes \xi (v_{i}^{\ast })
\end{equation*}%
est un isomorphisme d'espaces vectoriels.

\item Lorsque $A=\mathbb{D=}$ $\left\{ a+\varepsilon b:a,b\in 
\mathbb{R}
,\varepsilon ^{2}=0\right\} $ est l'ensemble des nombres duaux , qui est
isomorphe \`{a} l'alg\`{e}bre des polyn\^{o}mes tronqu\'{e}s $\mathbb{R}%
\left[ x\right] /\left( x^{2}\right) $, on note\ $\left( 1^{\ast
},\varepsilon ^{\ast }\right) $ la base duale de la base canonique $\left(
1\ ,\varepsilon \right) $\ de $\mathbb{D}$. La vari\'{e}t\'{e} $M^{\mathbb{D}%
}$ est identifi\'{e}e au fibt\'{e} tangent $TM$ par l'application 
\begin{equation*}
M^{\mathbb{D}}\longrightarrow TM,\xi \longmapsto \ \varepsilon ^{\ast }\circ
\xi \text{ }
\end{equation*}

L'application r\'{e}ciproque \'{e}tant 
\begin{equation*}
\ TM\longrightarrow M^{\mathbb{D}},v\longmapsto \left\{ \xi :f\longmapsto
f(p)+\varepsilon \cdot v(p)\right\} \ 
\end{equation*}%
si $v\in T_{p}M.$

\item Plus g\'{e}n\'{e}ralement, si $A=\mathbb{R}\left[ x_{1},...,x_{s}%
\right] /\left( x_{1},...,x_{s}\right) ^{k+1}$, alors $M^{A}=J_{0}^{k}(%
\mathbb{R}
^{s},M)$.
\end{enumerate}

\subsection{Champs de vecteurs sur $M^{A}$}

L'ensemble, $C^{\infty }(M^{A},A)$, des fonctions de classe $C^{\infty }$
sur $M^{A}$ \`{a} valeurs dans $A$ est une $A$-alg\`{e}bre commutative
unitaire, d'unit\'{e} l'application 
\begin{equation*}
1_{C^{\infty }(M^{A},A)}:\xi \longmapsto 1_{A}.
\end{equation*}

Pour $f\in C^{\infty }(M)$, l'application%
\begin{equation*}
f^{A}:M^{A}\longrightarrow A,\xi \longmapsto \xi (f)
\end{equation*}%
est de classe $C^{\infty }$ et l'application 
\begin{equation*}
\ C^{\infty }(M)\longrightarrow C^{\infty }(M^{A},A),f\longmapsto f^{A}
\end{equation*}%
est un homomorphisme d'alg\`{e}bres.

Si $(a_{i})_{i=1,...,s}$ est une base de $A$ et $(a_{i}^{\ast })_{i=1,...,s}$
la base duale, on identifie $C^{\infty }(M^{A},A)$ \`{a} $A\otimes C^{\infty
}(M^{A})$ par l'application 
\begin{equation*}
\sigma :\varphi \ \longmapsto \overset{s}{\underset{i=1}{\sum }}a_{i}\otimes
(a_{i}^{\ast }\circ \varphi )\text{.}
\end{equation*}%
Ainsi, $\ \sigma (f^{A})=\overset{s}{\underset{i=1}{\sum }}a_{i}\otimes
(a_{i}^{\ast }\circ f^{A})$ pour tout $f\in C^{\infty }(M)$.

On note 
\begin{equation*}
\gamma :C^{\infty }(M)\longrightarrow A\otimes C^{\infty
}(M^{A}),f\longmapsto \overset{s}{\underset{i=1}{\sum }}a_{i}\otimes
(a_{i}^{\ast }\circ f^{A})
\end{equation*}
et $Der_{\gamma }\left[ C^{\infty }(M),A\otimes C^{\infty }(M^{A})\right] $
le $A\otimes C^{\infty }(M^{A})$-module des $\gamma $-d\'{e}rivations de $%
C^{\infty }(M)$ dans $A\otimes C^{\infty }(M^{A})$ c'est-\`{a}-dire
l'ensemble des applications $%
\mathbb{R}
$-lin\'{e}aires 
\begin{equation*}
D:C^{\infty }(M)\longrightarrow A\otimes C^{\infty }(M^{A})
\end{equation*}%
telles que, pour $f$ et $g$ appartenant \`{a} $C^{\infty }(M)$,\qquad\ 
\begin{equation*}
D(fg)=D(f)\cdot \gamma (g)+\gamma (f)\cdot D(g)\text{.}
\end{equation*}

Une d\'{e}rivation de $C^{\infty }(M)$ dans $C^{\infty }(M^{A},A)$\ \cite%
{b-o}, est une d\'{e}rivation par rapport \`{a} l'homomorphisme 
\begin{equation*}
C^{\infty }(M)\longrightarrow C^{\infty }(M^{A},A),f\longmapsto f^{A}\text{ }
\end{equation*}%
c'est-\`{a}-dire, une application $%
\mathbb{R}
$- lin\'{e}aire 
\begin{equation*}
X:C^{\infty }(M)\longrightarrow C^{\infty }(M^{A},A)
\end{equation*}%
telle que, pour $f$ et $g$ appartenant \`{a} $C^{\infty }(M)$,%
\begin{equation*}
X(fg)=X(f)\cdot g^{A}+f^{A}\cdot X(g)\text{.}
\end{equation*}%
L'ensemble $Der\left[ C^{\infty }(M),C^{\infty }(M^{A},A)\right] $, des d%
\'{e}rivations de $C^{\infty }(M)$ dans $C^{\infty }(M^{A},A)$, est un $%
C^{\infty }(M^{A},A)$-module. D'apr\`{e}s \cite{ok1}, \cite{ok2},
l'application 
\begin{equation*}
Der\left[ C^{\infty }(M^{A})\right] \longrightarrow Der_{\gamma }\left[
C^{\infty }(M),A\otimes C^{\infty }(M^{A})\right] ,X\longmapsto
(id_{A}\otimes X)\circ \gamma ,
\end{equation*}%
est un isomorphisme de $C^{\infty }(M^{A})$-modules .

Il s'ensuit que l'application%
\begin{equation*}
Der\left[ C^{\infty }(M^{A})\right] \longrightarrow Der\left[ C^{\infty
}(M),C^{\infty }(M^{A},A)\right] ,X\longmapsto \sigma ^{-1}\circ
(id_{A}\otimes X)\circ \gamma ,
\end{equation*}%
est un isomorphisme de $C^{\infty }(M^{A})$-modules qui permet de
transporter sur $Der\left[ C^{\infty }(M^{A})\right] $ la structure de $%
C^{\infty }(M^{A},A)$-module de $Der\left[ C^{\infty }(M),C^{\infty
}(M^{A},A)\right] $. On peut donc regarder un champ de vecteurs sur $M^{A}$
comme une d\'{e}rivation de $C^{\infty }(M)$ dans $C^{\infty }(M^{A},A)$ 
\cite{b-o}

\begin{proposition}
Les assertions suivantes sont \'{e}quivalentes:

\begin{enumerate}
\item Un champ de vecteurs sur $M^{A}$ est une section diff\'{e}rentiable du
fibr\'{e} tangent $(TM^{A},\pi _{M^{A}},M^{A})$;

\item Un champ de vecteurs sur $M^{A}$ est une d\'{e}rivation de $C^{\infty
}(M^{A})$;

\item Un champ de vecteurs sur $M^{A}$ est une d\'{e}rivation de $C^{\infty
}(M)$ dans $C^{\infty }(M^{A},A)$.
\end{enumerate}
\end{proposition}

\begin{example}
\bigskip Soit $C$ est le champ de Liouville sur le fibr\'{e} tangent $TM$.
Dans un syst\`{e}me de coordonn\'{e}es locales $(x_{1},...,x_{n})$ sur la
vari\'{e}t\'{e} $M$, si $y_{i}=dx_{i}$ d\'{e}signe la coordonn\'{e}e sur la
fibre, 
\begin{equation*}
C=\underset{i=1}{\overset{n}{\dsum }}\text{ }y_{i}\frac{\partial }{\partial
y_{i}}
\end{equation*}%
cest-\`{a}-dire $C(x_{i})=0$ et $\ C(y_{i})=\ y_{i}$.

Puisque l'application 
\begin{equation*}
\gamma :C^{\infty }(M)\longrightarrow \mathbb{D}\otimes C^{\infty
}(TM),f\longmapsto 1\otimes 1^{\ast }\circ f^{\mathbb{D}}+\varepsilon
\otimes \varepsilon ^{\ast }\circ f^{\mathbb{D}}
\end{equation*}%
o\`{u} 
\begin{equation*}
f^{\mathbb{D}}\left( v\right) =f(p)+\varepsilon \cdot v(f)\text{, }v\in TpM
\end{equation*}%
est telle que 
\begin{eqnarray*}
\gamma \left( x_{i}\right) \left( 1\otimes v\right) &=&\left[ 1\otimes
1^{\ast }\circ x_{i}^{\mathbb{D}}+\varepsilon \otimes \varepsilon ^{\ast
}\circ x_{i}^{\mathbb{D}}\right] \left( \left( 1\otimes v\right) \right) \\
&=&1\otimes 1^{\ast }\left[ x_{i}(p)+\varepsilon \cdot v(x_{i})\right]
+\varepsilon \otimes \varepsilon ^{\ast }\left[ x_{i}(p)+\varepsilon \cdot
v(x_{i})\right] \\
&=&1\otimes \ x_{i}(p)+\varepsilon \otimes v(x_{i}) \\
&=&\left[ 1\otimes x_{i}+\varepsilon \otimes dx_{i}\right] \left( \left(
1\otimes v\right) \right) \\
&=&\left[ 1\otimes x_{i}+\varepsilon \otimes y_{i}\right] \left( 1\otimes
v\right) \text{,}
\end{eqnarray*}%
c'est-\`{a}-dire 
\begin{equation*}
\gamma \left( x_{i}\right) =1\otimes x_{i}+\varepsilon \otimes y_{i}\text{,}
\end{equation*}%
la d\'{e}rivation $X:C^{\infty }(M)\longrightarrow C^{\infty }(TM,\mathbb{D}%
) $ correspondant au champ de Liouville $C$ est donn\'{e}e par%
\begin{eqnarray*}
X\left( x_{i}\right) &=&\sigma ^{-1}\circ \left( id_{\mathbb{D}}\otimes
C\right) \circ \gamma \left( x_{i}\right) \\
&=&\sigma ^{-1}\circ \left( id_{\mathbb{D}}\otimes C\right) \circ \left(
1\otimes x_{i}+\varepsilon \otimes y_{i}\right) \\
&=&\sigma ^{-1}\circ \left[ 1\otimes C(x_{i})+\varepsilon \otimes C(y_{i})%
\right] \\
&=&\ \varepsilon \cdot y_{i}\text{.}
\end{eqnarray*}
\end{example}

Dans toute la suite, nous regarderons un champ de vecteurs comme une d\'{e}%
rivation de $C^{\infty }(M)$ dans $C^{\infty }(M^{A},A)$.

\subsubsection{\textbf{\ Champs de vecteurs sur }$M^{A}$\textbf{\ provenant
des d\'{e}rivations de }$A$}

\begin{proposition}
Si $d$ est une d\'{e}rivation de $A$, alors l'application \ 
\begin{equation*}
d^{\ast }:C^{\infty }(M)\longrightarrow C^{\infty }(M^{A},A),f\longmapsto
(-d)\circ f^{A},
\end{equation*}%
est un champ de vecteurs sur $M^{A}$.
\end{proposition}

\begin{demonstration}
L'application $d^{\ast }$ est $%
\mathbb{R}
$-lin\'{e}aire et pour $f$ et $g$ appartenant \`{a} $C^{\infty }(M)$ et pour 
$\xi \in M^{A}$, on a:%
\begin{eqnarray*}
d^{\ast }(fg)(\xi ) &=&(-d)\circ (fg)^{A}(\xi ) \\
&=&(-d)\circ (f^{A}\cdot g^{A})(\xi ) \\
&=&(-d)\left[ f^{A}(\xi )\cdot g^{A}(\xi )\right] \  \\
&=&(-d)\left[ f^{A}(\xi )\right] \cdot g^{A}(\xi )+f^{A}(\xi )\cdot (-d)%
\left[ g^{A}(\xi )\right] \\
&=&\left[ (-d)\circ f^{A}\right] (\xi )\cdot g^{A}(\xi )+f^{A}(\xi )\cdot %
\left[ (-d)\circ f^{A}\right] (\xi ) \\
&=&(\left[ (-d)\circ f^{A}\right] \cdot g^{A}+f^{A}\cdot \left[ (-d)\circ
f^{A}\right] )(\xi ) \\
&=&\left[ d^{\ast }(f)\cdot g^{A}+f^{A}\cdot d^{\ast }(g)\right] (\xi )\text{%
.}
\end{eqnarray*}%
Comme $\xi $ est quelconque, on d\'{e}duit que 
\begin{equation*}
d^{\ast }(fg)=d^{\ast }(f)\cdot g^{A}+f^{A}\cdot d^{\ast }(g)\text{.}
\end{equation*}%
Ainsi, $d^{\ast }$ est un champ de vecteurs sur $M^{A}$.
\end{demonstration}

Le champ de vecteurs $d^{\ast }$ est le champ de vecteurs sur $M^{A}$ associ%
\'{e} \`{a} la d\'{e}rivation $d$ de $A$ et on a pour $d_{1}$, $d_{2}$ , $d$
trois d\'{e}rivations de $A$ et pour $a$ $\in A$, on a \cite{b-o} :%
\begin{align*}
\left[ d_{1}^{\ast },d_{2}^{\ast }\right] & =\left[ d_{1},d_{2}\right]
^{\ast }\text{;} \\
(a\cdot d)^{\ast }& =a\cdot d^{\ast }\text{.}
\end{align*}

\section{Feuilletage induit par les d\'{e}rivations de $A$}

\begin{theoreme}
Soient $d_{1},\ ...,d_{r}$ une base de $Der(A)$. Les champs de vecteurs, $%
d_{1}^{\ast },\ ...,d_{r}^{\ast }$ induits par $d_{1},\ ...,d_{r}$ sur $%
M^{A} $ sont lin\'{e}airement ind\'{e}pendants et d\'{e}finissent un
feuilletage $\mathcal{F}_{r}$ de dimension $r$ sur $M^{A}$.
\end{theoreme}

\begin{demonstration}
Soient les applications $\varphi _{1},\ ...,\varphi _{r}$ $\in C^{\infty
}(M^{A},A)$ telles que $\varphi _{1}\cdot d_{1}^{\ast }+...+\varphi
_{r}\cdot d_{r}^{\ast }=0$.

Pour $f\in C^{\infty }(M)$, et pour tout $\xi \in M^{A}$, on a 
\begin{align*}
0& =\left[ \varphi _{1}\cdot d_{1}^{\ast }\ +...+\varphi _{r}\cdot
d_{r}^{\ast }\right] (f)\left( \xi \right) \\
& =\left[ \varphi _{1}\cdot \left( (-d_{1})\circ f^{A}\right) +...+\varphi
_{r}\cdot \left( (-d_{r})\circ f^{A}\right) \right] \left( \xi \right) \\
& =-\varphi _{1}\left( \xi \right) \cdot d_{1}\left( \xi \left( f\right)
\right) -...-\varphi _{r}\left( \xi \right) \cdot d_{r}\left( \xi \left(
f\right) \right) \\
& =-\left[ \varphi _{1}\left( \xi \right) \cdot d_{1}\left( \xi \left(
f\right) \right) +...+\varphi _{r}\left( \xi \right) \cdot d_{r}\left( \xi
\left( f\right) \right) \right]
\end{align*}%
Puisque les d\'{e}rivations $d_{1},\ ...,d_{r}$ sont lin\'{e}airements ind%
\'{e}pendantes, on a alors 
\begin{equation*}
\varphi _{1}\left( \xi \right) =...=\varphi _{r}\left( \xi \right) =0\text{.}
\end{equation*}%
Comme $\xi $ est quelconque, on conclut que 
\begin{equation*}
\varphi _{1}\ =...=\varphi _{r}\ =0
\end{equation*}%
c'est-\`{a}-dire que les champs de vecteurs $d_{1}^{\ast },\ ...,d_{r}^{\ast
}$ sont lin\'{e}airement ind\'{e}pendantes.

De plus, pour $i,j\in \left\{ 1,...,r\right\} $, 
\begin{equation*}
\left[ d_{i}^{\ast },d_{j}^{\ast }\right] =\left[ d_{i},d_{j}\right] ^{\ast }%
\text{.}
\end{equation*}%
Le syst\`{e}me diff\'{e}rentiel engendr\'{e} par $d_{1}^{\ast },\
...,d_{r}^{\ast }$ est donc compl\'{e}tement int\'{e}grable. Il d\'{e}finit
donc un feuilletage $\mathcal{F}_{r}$ de dimension $r$ que nous appelons
feuilletage canonique sur $M^{A}$.

Ce qui ach\`{e}ve la d\'{e}monstration.
\end{demonstration}

\subsection{Feuilletage canonique sur le fibr\'{e} tangent $TM$}

\begin{proposition}
Le feuilletage canonique sur le fibr\'{e} tangent $TM$ est le feuilletage d%
\'{e}fini par son champ de vecteurs canonique (champ de Liouville) $C$.
\end{proposition}

\begin{demonstration}
Soit $d$ une d\'{e}rivation de $\mathbb{D}$, alors 
\begin{eqnarray*}
0 &=&d\left( \varepsilon ^{2}\right) \\
&=&d\left( \varepsilon \cdot \varepsilon \right) \\
&=&2\varepsilon \cdot d\left( \varepsilon \right)
\end{eqnarray*}%
Il existe donc $\lambda \in 
\mathbb{R}
$ tel que $d\left( \varepsilon \right) =\lambda \varepsilon $. Toute d\'{e}%
rivation de $\mathbb{D}$ est donc de la forme telle que $d(\varepsilon
)=\lambda \varepsilon $ o\`{u} $\lambda \in 
\mathbb{R}
$.

Soit $d_{0}$ la d\'{e}rivation de $\mathbb{D}$ telle que $d_{0}\left(
\varepsilon \right) =-\varepsilon $, alors $Der(\mathbb{D)=%
\mathbb{R}
\cdot }d_{0}$ est l'alg\`{e}bre de Lie de dimension $1$, engendr\'{e}e par $%
d_{0}$.

Le champ de vecteurs provenant de la d\'{e}rivation $d_{0}$ est l'application%
$\ $%
\begin{equation*}
d_{0}^{\ast }:C^{\infty }(M)\longrightarrow C^{\infty }(TM,\mathbb{D}%
),f\longmapsto (-d_{0})\circ f^{\mathbb{D}},
\end{equation*}%
c'est-\`{a}-dire, si $(x_{1},...,x_{n})$ est un syst\`{e}me de coordonn\'{e}%
es locales sur la vari\'{e}t\'{e} $M$, et si $(y_{1},...,y_{n})$ d\'{e}%
signent les coordonn\'{e}es sur la fibre, 
\begin{eqnarray*}
(d_{0}^{\ast })(x_{i}) &=&(-d_{0})\circ \left( x_{i}\right) ^{\mathbb{D}} \\
&=&-d_{0}\circ \left( x_{i}+\varepsilon y_{i}\right) \\
&=&-d_{0}\circ \left[ \varepsilon \cdot y_{i})\right] \\
\ &=&\varepsilon \cdot y_{i}
\end{eqnarray*}%
Le champ de vecteurs $d_{0}^{\ast }$ est donc le champ de Liouville sur le
fibr\'{e} tangent et le feuilletage induit par $d_{0}^{\ast }$ sur $TM$ est
\ par cons\'{e}quent le feuilletage canonique de $TM$.
\end{demonstration}

\

\end{document}